\documentclass[dvips,noinfoline,ejs]{imsart}

\RequirePackage[OT1]{fontenc}
\RequirePackage{amsthm,amsmath}
\usepackage[square]{natbib}
\RequirePackage[colorlinks]{hyperref}
\RequirePackage{hypernat}

\doi{10.1214/08-EJS228}
\pubyear{2008}
\volume{2}
\firstpage{542}
\lastpage{563}

\def\1{{\mathbbm 1}}
\def\N{{\mathbb{N}}}

\def\R{{\mathbb{R}}}

\def\E{{\mathbb{E}}}

\def\var{{\mathrm{Var}}}
\def\pen{{\mathrm{pen}}}
\def\eps{{\varepsilon}}

\def\M{\mathcal{M}}
\def\G{\mathcal{G}}

\def\p{\mathbb{P}}

\def\crit{{\mathrm{Crit}}}

\def\EH{{\mathrm{\sf EDkhi}}}
\def\H{\mathrm{\sf Dkhi}}

\def\NN{\mathcal{N}}

\def\x{\mathrm{X}}
\def\g{{\bf g}}
\def\thetah{\hat\theta}

\newcommand{\eref}[1]{(\ref{#1})}
\newcommand{\pa}[1]{\left({#1}\right)}
\newcommand{\pab}[1]{\big({#1}\big)}
\newcommand{\paB}[1]{\Big({#1}\Big)}

\newcommand{\crob}[1]{\big[{#1}\big]}
\newcommand{\cro}[1]{\left[{#1}\right]}

\newcommand{\acb}[1]{\big\{{#1}\big\}}
\newcommand{\ac}[1]{\left\{{#1}\right\}}

\RequirePackage{amssymb}

\newtheorem{lemma}{Lemma}
\newtheorem{theorem}{Theorem}
\newtheorem{proposition}{Proposition}
\newtheorem{corollary}{Corollary}

\RequirePackage{amssymb}

\begin{document}

\begin{frontmatter}
\title{Estimation of Gaussian graphs by model selection}
\runtitle{Estimation of Gaussian graphs}

\begin{aug}
\author{\fnms{Christophe} \snm{Giraud}\ead[label=e1]{cgiraud@math.unice.fr}}
\address{Universit\'e de Nice, Laboratoire J.A.\;Dieudonn\'e,  06108 Nice cedex 02, France  \\
INRA, MIA 78352 Jouy-en-Josas Cedex, France\\ \printead{e1}}
\runauthor{C. Giraud}
\end{aug}

\begin{abstract}
We investigate in this paper the estimation of  Gaussian graphs by
model selection from a non-asymptotic point of view. We start from an
$n$-sample of a Gaussian law $\p_{C}$ in $\R^p$ and  focus on the
disadvantageous case where $n$ is smaller than $p$. To estimate the
graph of conditional dependences of $\p_{C}$, we introduce a collection
of candidate graphs and then select  one of them by minimizing a
penalized empirical risk. Our main result assesses the performance of
the procedure in a non-asymptotic setting. We pay special attention to
the maximal degree $D$ of the graphs that we can handle, which turns to
be roughly $n/(2\log p)$.
\end{abstract}



\received{\smonth{4} \syear{2008}}

\end{frontmatter}


\section{Introduction}

Let us consider a Gaussian law $\p_{C}$ in $\R^p$ with mean 0 and
positive definite covariance matrix $C$. We write $\theta$ for the
matrix of the  regression coefficients  associated to the law $\p_{C}$,
more precisely  $\theta=\crob{\theta^{(j)}_{i}}_{i,j=1,\ldots,p}$ is the
$p\times p$ matrix such that $\theta_{j}^{(j)}=0$ for $j=1,\ldots,p$
and
$$\E\crob{\x^{(j)}\,\big|\,\x^{(k)},\ k\neq j}=\sum_{k\neq j}\theta_{k}^{(j)}\x^{(k)},\quad j\in \ac{1,\ldots,p},\quad \textrm{a.s.}$$
for any  random vector $\x=\pa{\x^{(1)},\ldots,\x^{(p)}}^T$ of law
$\p_{C}$. Our aim is to estimate the matrix $\theta$  by model
selection from an $n$-sample $X_{1},\ldots,X_{n}$ i.i.d. with law
$\p_{C}$. We will focus on the disadvantageous case where the sample
size $n$ is smaller than the dimension $p$.

We call henceforth \emph{shape of} $\theta$, the set of the couples  of
integers $(i,j)$ such that $\theta_{i}^{(j)}\neq 0$. The shape of
$\theta$ is usually represented by a graph $\g$ with $p$ labeled
vertices $\ac{1,\ldots,p}$, by  setting an edge between the vertices
$i$ and $j$ when $\theta_{i}^{(j)}\neq 0$. This graph is well-defined
since $\theta_{i}^{(j)}=0$ if and only if $\theta_{j}^{(i)}=0$; the
latter property may be seen e.g.\;on the formula
$\theta_{i}^{(j)}=-(C^{-1})_{i,j}/(C^{-1})_{j,j}$ for all $i\neq j$.
The graph $\g$ is of interest for the statistician since it depicts the
conditional dependences of the variables $X^{(j)}$s. Actually, there is
an edge between $i$ and $j$ if and only if $X^{(i)}$ \emph{is not
independent of $X^{(j)}$ conditionally on the other variables.} The
objective in Gaussian graphs estimation is usually to detect the graph
$\g$. Even if the purpose of our procedure is to estimate $\theta$ and
not $\g$, we propose to simultaneously estimate $\g$  as follows. We
associate with our estimator $\hat \theta$ of $\theta$, the graph $\hat
\g$ where we set an edge between the vertices $i$ and $j$ when
$\thetah_{i}^{(j)}$  is non-zero.

Estimation of Gaussian graphs with $n\ll p$ is a current active field
of research motivated by applications in postgenomic. Biotechnological
developments (microarrays, 2D-electrophoresis, etc) enable to produce a
huge amount of proteomic and transcriptomic data. One of the challenge
in postgenomic  is to infer from these data the regulation network of a
family of  genes (or proteins). The task is challenging for the
statistician due to the very high-dimensional nature of the data and
the small sample size. For example, microarrays measure the expression
levels of a few thousand genes (typically 4000) and the sample size $n$
is no more than a few tens. The Gaussian graphical modeling appears to
be a valuable tool for this issue, see the papers of
Kishino and  Waddell~\cite{KW}, Dobra {\it et al}~\cite{Detal}, Wu and
Ye~\cite{WY}.  The gene expression levels in the microarray are modeled
by a  Gaussian law $\p_{C}$ and   the regulation network of the genes
is then depicted by the graph $\g$ of the conditional dependences.

Various procedures have been proposed to perform graph estimation when
$p>n$. Many are based on multiple testing, see for instance the papers
of Sch\"afer and Strimmer~\cite{SS}, Drton and Perlman~\cite{DP04,
DP07} or Wille and B\"uhlmann~\cite{WB}. We also mention the work of
Verzelen and Villers~\cite{VV} for testing in a non-asymptotic
framework whether  there are (or not) missing edges in a given graph.
Recently, several authors  advocate to take advantage of the nice
computational properties of the $l^1$-penalization to either estimate
the graph $\g$ or the concentration matrix $C^{-1}$. Meinshausen and
B\"uhlmann~\cite{MB} propose to learn the graph $\g$ by regressing with
the Lasso  each variable against the others. Huang {\it et
al.}~\cite{HLPL} or Yuan and Lin~\cite{YL}  (see also Banerjee {\it et
al.}~\cite{BGA} and Friedman {\it et al.}~\cite{FHT}) suggest in turn
to rather estimate $C^{-1}$ by  minimizing the log-likelihood for the
concentration matrix penalized by the $l^1$-norm. The performance of
these algorithms are mostly unknown: the  few theoretical results are
only valid  under restrictive conditions on the covariance matrix and
for large $n$ (asymptotic setting). In addition to these few
theoretical results, Villers {\it et al.}  \cite{HSV} propose a
numerical investigation of the validity domain of some of the above
mentioned procedures.

Our aim in this work is to investigate Gaussian graph estimation  by
model selection from a non-asymptotic point of view. We propose a
procedure to estimate $\theta$ and  assess its performance in a
non-asymptotic setting.
Then, we discuss on the maximum degree of the graphs that we can
accurately estimate and explore the performance of our estimation
procedure  in a small numerical study.

We will use the Mean Square Error of Prediction (MSEP) as a criterion
to assess the quality of our procedure. To define this quantity, we
introduce a few notations. For any $k,q\in\N$, we  write
$\|\cdot\|_{k\times q}$ for the Frobenius  norm in $\R^{k\times q}$,
namely $\|A\|^2_{k\times q}={\rm Trace}\,(A^TA)$,
for any $A\in\R^{k\times q}.$
The MSEP of the estimator $\thetah$ is then
$$\textrm{MSEP}(\thetah)=\E\cro{\|C^{1/2}(\thetah-\theta)\|_{p\times p}^2}=\E\cro{\|\x_{new}^T(\thetah-\theta)\|_{1\times p}^2},$$
where $C^{1/2}$ is the positive square root of $C$ and $\x_{new}$ is a
random vector, independent of $\thetah$, with distribution $\p_{C}$. We
underline that the MSEP focus on the quality of the estimation of
$\theta$ and not of $\g$. In particular, we do not aim to estimate at
best the ``true'' graph $\g$, but rather to estimate at best the
regression matrix $\theta$. We choose this point of view for two
reasons. First, we do not believe that the matrix $\theta$ is exactly
sparse in practice, in the sense that $\theta_{i}^{(j)}=0$ for most of
the $i,j\in\ac{1,\ldots,p}$. Rather, we want to handle cases where the
matrix $\theta$ is only approximately sparse, which means that there
exists a sparse matrix $\theta^*$ which is a good approximation of
$\theta$. In this case, the shape $\g$ of $\theta$ may not be sparse at
all, it can even be the complete graph. Our goal is then not to
estimate $\g$ but rather to capture the main conditional dependences
given by the shape $\g^*$ of $\theta^*$. The second reason for
considering the MSEP as a quality criterion for our procedure is that
we want to quantify the fact that we do not want to miss the important
conditional dependences, but we do not worry too much missing a weak
one. In other words, even in the case where the shape $\g$ of $\theta$
is sparse, we are interested in finding the main edges of $\g$
(corresponding to strong conditional dependences) and we do not really
care of missing a ``weak'' edge which is overwhelmed by the noise. The
MSEP is a possible way to take this issue into account.

To estimate $\theta$, we will first introduce a collection $\M$ of
graphs, which are our candidates for describing the shape $\g$ of
$\theta$. If we have no prior information on $\g$, a possible choice
for $\M$ is the set of all graphs with degree\footnote{the degree of a
graph corresponds to the maximum number of edges incident to a vertex.}
less than some fixed integer $D$. Then, we associate with each graph
$m\in\M$, an estimator $\thetah_{m}$ of $\theta$ by minimizing an
empirical version of the MSEP with the constraint that the shape of
$\thetah_{m}$ is given by $m$, see Section~\ref{procedure} for the
details. Finally, we select one of the candidate graph $\hat m$ by
minimizing a penalized empirical MSEP and set $\thetah=\thetah_{\hat
m}$. Our main result roughly states that when the candidate graphs have
a degree smaller than  $n/(2\log p)$,  the MSEP of $\hat \theta$ nearly
achieves, up to a $\log (p)$ factor, the minimal MSEP of the collection
of estimators $\{\thetah_{m},\ m\in\M\}$.

It is of practical interest to know if the  condition on the degree of
the candidate graphs  can be avoided. This point is discussed in
Section~\ref{optimal}, where we emphasize  that it is hopeless to try
to estimate accurately graphs with a degree $D$ large compared to
$n/(1+\log (p/n))$. We also prove that the size of the penalty involved
in the selection procedure  is minimal in some sense.

The remaining of the paper is organized as follows. After introducing a
few notations, we describe the estimation procedure  in
Section~\ref{procedure} and  state our main results  in
Section~\ref{perf}. Section~\ref{numerique} is devoted to a small
numerical study and Section~\ref{proofs} to the proofs.

\subsection*{A few notations}
Before describing our estimation procedure, we introduce a few notations about graphs we shall use all along the paper.

\subsubsection*{a. Graphs}

The set of the graphs with $p$ vertices labeled by $\ac{1,\ldots,p}$ is
in bijection  with the set $\G$ of all the subset $g$ of
$\ac{1,\ldots,p}^2$ fulfilling
\begin{itemize}
\item $(j,j)\notin g$ for all $j\in\ac{1,\ldots,p}$,
\item $(i,j)\in g\ \Rightarrow (j,i)\in g$ for all $i,j\in\ac{1,\ldots,p}$.
\end{itemize}
Indeed, to any $g\in\G$ we can associate a graph with $p$ vertices
labeled by $\ac{1,\ldots,p}$ by setting an edge between the vertices
$i$ and $j$ if and only if $(i,j)\in g$. For simplicity, we call
henceforth ``graph'' any element $g$ of $\G$.

For a graph $g\in\G$ and an integer $j\in\ac{1,\ldots,p}$, we set
$g_{j}=\ac{i:\ (i,j)\in g}$ and denote by $|g_{j}|$ the cardinality of
$g_{j}$. Finally, we define the degree of $g$ by ${\rm
deg}(g)=\max\ac{|g_{j}|:\ j=1,\ldots,p}$.

\subsubsection*{b. Directed graphs}

As before, we will represent the set of the directed graph with $p$
vertices labeled by $\ac{1,\ldots,p}$ by the set $\G^+$ of all the
subset $g$ of $\ac{1,\ldots,p}^2$ fulfilling $(j,j)\notin g$ for all
$j\in\ac{1,\ldots,p}$. More precisely, we associate with  $g\in\G^+$
the directed graph with $p$ vertices labeled by $\ac{1,\ldots,p}$ and
with directed edges from $i$ to $j$ if and only if $(i,j)\in g$.

We note that $\G\subset \G^+$ and we extend to $g\in\G^+$ the above
definitions of $g_{j}$, $|g_{j}|$, and ${\rm deg}(g)$. Although $\G$ is
contained in $\G^+$, it should be noted that the associated
interpretation is different since the graphs in $\G^+$ are directed
with possibly two directed edges between two vertices.

\section{Estimation procedure}\label{procedure}

In this section, we explain our procedure to estimate $\theta$. We
first introduce a collection of graphs and models, then we associate
with each model an estimator and finally we give a  procedure to select
one of them.

\subsection{Collection of graphs and models}\label{candidates}

Our estimation procedure starts with the choice of either a collection
$\M\subset\G$ of graphs or a collection $\M\subset\G^+$  of directed
graphs which are our candidates to describe the shape of $\theta$.
Among the possible choices for $\M$ we mention four of them:
\begin{enumerate}
\item the set $\M^{\#}_{D}\subset\G$ of all  graphs with at most $D$ edges,
\item the set $\M^{{\rm deg}}_{D}\subset \G$ of all  graphs with degree less than $D$,
\eject
\item the set $\M^{\#,+}_{D}\subset\G^+$ of all  directed graphs with at most $D$ directed edges,
\item the set $\M^{{\rm deg},+}_{D}\subset \G^+$ of all  directed graphs with degree less than $D$.
\end{enumerate}
We call degree of $\M$ the integer $D_{\M}=\max\ac{\deg(m)\, m\in\M}$ and
note that the above collections of graphs have a degree bounded by $D$.

To the collection of  graphs $\M$, we associate the following
collection  $\{\Theta_{m},\break m\in\M\}$ of models  to estimate $\theta$.
The model $\Theta_{m}$ is the linear space of those matrices in
$\R^{p\times p}$ whose shape is given by the graph $m$, namely
$$\Theta_{m}=\acb{A\in\R^{p\times p}:\ (i,j)\notin m \Rightarrow A_{i}^{(j)}=0}.$$

As mentioned before, we known that $\theta_{i}^{(j)}=0$ if and only if
$\theta_{j}^{(i)}=0$, so it seems irrelevant to (possibly) introduce
directed graphs instead of graphs. Nevertheless, we must keep in mind
that our aim is to estimate $\theta$ at best in terms of the MSEP. In
some cases, the results can be improved when using directed graphs
instead of graphs, typically when for some $i,j\in\ac{1,\ldots,p}$ the
variance of $\theta_{i}^{(j)}\x^{(i)}$ is large compared to the
conditional variance $\var(\x^{(j)}| \x^{(k)},\ k\neq j)$, whereas the
variance of $\theta^{(i)}_{j}\x^{(j)}$ is small compared to
$\var(\x^{(i)}| \x^{(k)},\ k\neq i)$. Finally, we note the following
inclusions for the families of models mentioned above
$$\bigcup_{m\in\M^{\#,+}_{D}}\Theta_{m}\ \subset\ \bigcup_{m\in\M^{\#}_{D}}\Theta_{m}\ \subset\ \bigcup_{m\in\M^{{\rm deg}}_{D}}\Theta_{m}\ \subset\  \bigcup_{m\in\M^{{\rm deg},+}_{D}}\Theta_{m}.$$

\subsection{Collection of estimators}
We assume henceforth that $3\leq n<p$ and that the degree $D_{\M}$ of
$\M$ is upper bounded by some integer $D\leq n-2$. We start with  $n$
observations $X_{1},\ldots,X_{n}$ i.i.d.\;with law $\p_{C}$ and we
denote by $X$ the $n\times p$ matrix $X=\cro{X_{1},\ldots,X_{n}}^T$. In
the following, we write $A^{(1)},\ldots,A^{(p)}$  for the $p$ columns
of a matrix $A\in\R^{k\times p}$.

We remind the reader that
$\|C^{1/2}(I-\theta)\|^2=\inf_{A\in\Theta}\|C^{1/2}(I-A)\|^2,$
where $\Theta$ is the space of $p\times p$ matrices with 0 on the diagonal. An empirical version of $\|C^{1/2}(I-A)\|^2$ is $n^{-1}\|X(I-A)\|^2_{n\times p}$, which can also be viewed as an empirical version of the loss $\|C^{1/2}(A-\theta)\|^2$,
since by Pythagorean theorem
$\|C^{1/2}(A-\theta)\|^2=\|C^{1/2}(I-A)\|^2-\|C^{1/2}(I-\theta)\|^2$, for all $A\in\Theta.$

In this direction, we associate with any $m\in \M$, an estimator $\thetah_{m}$ of $\theta$ by minimizing on $\Theta_{m}$ this empirical risk
\begin{equation}\label{estimateur}
\|X(I-\thetah_{m})\|^2_{n\times p}=\min_{A\in \Theta_{m}}\|X(I-A)\|^2_{n \times p}.
\end{equation}
We note that the $p\times p$ matrix $\thetah_{m}$ then fulfills the equalities
$$X\thetah_{m}^{(j)}=\textrm{Proj}_{X\Theta^{(j)}_{m}}\pab{X^{(j)}},\quad\textrm{for } j=1,\ldots, p,$$
where $\Theta^{(j)}_{m}$ is the linear space $\Theta^{(j)}_{m}=\ac{\theta^{(j)}:\ \theta\in\Theta_{m}}\subset \R^p$ and $\textrm{Proj}_{X\Theta_{m}^{(j)}}$ is the orthogonal projector onto $X\Theta_{m}^{(j)}$ in $\R^n$ (for the usual scalar product). Hence,  since the covariance matrix $C$ is positive definite and $D$ is less than $n$, the minimizer of \eref{estimateur} is unique a.s.

\subsection{Selection procedure}
To estimate $\theta$, we will select one of the estimator $\thetah_{m}$ by minimizing some penalized version of the empirical risk $\|X(I-\thetah_{m})\|^2/n$. More precisely, we set $\thetah=\thetah_{\hat m}$ where $\hat m$ is any minimizer  on $\M$ of the criterion
\begin{equation}\label{selection}
\crit(m)=\sum_{j=1}^p\cro{\|X^{(j)}-X\thetah_{m}^{(j)}\|^2\times \pa{1+{\pen(|m_{j}|)    \over n-|m_{j}|}}   },
\end{equation}
with the penalty function $\pen:\N\to\R^+$ of the form of the penalties introduced   in Baraud {\it et al.}~\cite{BGH}. To compute this penalty, we define for any integers $d$ and $N$ the $\H$ function by
$$\H(d,N,x)=
\p\pa{F_{d+2,N}\geq {x \over d+2}}-{x \over d}\,\p\pa{F_{d,N+2}\geq {N+2\over Nd}\, x},\quad x>0,$$
where $F_{d,N}$ denotes a Fisher random variable with $d$ and $N$ degrees of freedom. The function $x\mapsto \H(d,N,x)$ is decreasing and we write $\EH[d,N,x]$ for its inverse, see~\cite{BGH} Section 6.1 for details. Then, we fix some constant $K>1$ and set
\begin{equation}\label{new}
\pen(d)=K\,{n-d\over n-d-1}\,\EH\cro{d+1,n-d-1,\pa{C_{p-1}^d(d+1)^2}^{-1}}.
\end{equation}

\subsubsection*{Size of the penalty}
The size of the penalty $\pen(d)$ is roughly
$ 2K d \log p $ for large values of $p$. Indeed, we will work  in the sequel with collections of models, such that
$$ D_{\M}\leq \eta \;{n\over 2\pa{1.1+\sqrt{\log p}}^2},\quad \textrm{for some }\eta<1,$$
and then, we approximately have  for large values of $p$ and $n$
$$\pen(d)\lesssim K \pa{1+e^{\eta}\sqrt{2\log p}}^2(d+1),\quad d\in\ac{0,\ldots,D_{\M}},$$
see Proposition~4 in Baraud {\it et al.}~\cite{BGH} for an exact bound. In Section~\ref{min}, we show that the size of this penalty is minimal in some sense.

\subsubsection*{Choice of the tuning parameter $K$}
Increasing the value of $K$ decreases the size of the graph $\hat m$ that is selected. The choice $K=2$ gives  good control of the MSEP of $\thetah$, both theoretically and numerically (see Section~\ref{perf} and~\ref{numerique}). If we want that the rate of false discovery of edges remains smaller than $5\%$, the choice $K=3$ may also be appropriated.

\subsubsection*{Computational cost}
The computational cost of the selection procedure appears to be very high. For example, if $\M=\M^{{\rm deg},+}_{D}$ the computational complexity of the procedure increases as $p^{(D+1)}$ with the dimension $p$. In  a future work~\cite{GHV}, we will propose a modified version of this procedure, which presents a much smaller complexity.

\subsubsection*{A few additional remarks on the estimation procedure}
{\bf 1-} The matrix $\theta$ belongs to the set
$$\Gamma = \ac{\theta \in \R^{p\times p}:\ \exists\ K \textrm{ positive definite such that } \,\theta_{i,j}=-K_{i,j}/K_{j,j},\textrm{ for } i \neq j },$$
but the estimator $\thetah$  has no reason to belong to this space. To avoid this unpleasant feature, it would be natural to minimize~\eref{estimateur}
on the space $\Theta_{m}\cap\Gamma$ instead of $\Theta_{m}$. Unfortunately, we do not know how to handle this case neither theoretically nor numerically. We also emphasize that the matrix $\theta$ is not assumed to be exactly sparse, so it does not belong to any of the $\ac{\Theta_{m}\cap\Gamma,\ m\in\M}$ in general. In particular, it is unclear whether the MSEP of the estimator obtained by minimizing~\eref{estimateur} on $\Theta_{m}\cap\Gamma$ is smaller than the MSEP of $\thetah_{m}$.

{\bf 2-} In the special case where $\M=\M^{{\rm deg},+}_{D}$, the minimization of~\eref{selection} can be obtained by minimizing $\|X^{(j)}-X\thetah_{m}^{(j)}\|^2\times \pa{1+{\pen(|m_{j}|)    \over n-|m_{j}|}}$  independently for each $j$. This nice computational feature does not hold for the other collections of  graphs introduced in Section~\ref{candidates}.

\section{The main result}\label{perf}
 Next theorem gives an upper-bound on the MSEP of  a slight variation $\tilde\theta$ of $\thetah$, defined by
\begin{equation}\label{estimateur2}
\tilde\theta^{(j)}=\thetah^{(j)}\,{\bf 1}_{\ac{\|\thetah^{(j)}\|\leq \sqrt{p}\,T_{n}}},\ \textrm{for all } j\in\ac{1,\ldots,p}, \quad\textrm{with }\ T_{n}=n^{2 \log n}.
\end{equation}
We note that $\thetah$ and $\tilde\theta$ coincide  in practice since the threshold level $T_{n}$ increases very fast with $n$, e.g.\;$T_{20}\approx 6. 10^7$.

In the sequel, we write $\sigma^2_{j}=\pa{C^{-1}_{j,j}}^{-1}=\var(\x^{(j)}\mid \x^{(k)}, k\neq j)$ and define $\theta_{m}$  by
$$\| C^{1/2}(\theta-\theta_{m})\|^2=\min_{A_{m}\in\Theta_{m}}\| C^{1/2}(\theta-A_{m})\|^2.$$
\begin{theorem}\label{main}
Assume that $p>n\geq 3$ and $D_{\M}=\max\ac{{\rm deg}(m),\, m\in\M}$ fulfills the condition
\begin{equation}\label{conditionD}
1\leq D_{\M}\leq \eta \;{n\over 2\pa{1.1+\sqrt{\log p}}^2},\quad \textrm{for some }\eta<1.
\end{equation}
Then,  the MSEP of the estimator $\tilde\theta$ defined by~\eref{estimateur2}  is upper bounded by
\begin{eqnarray}\label{oracle1}
\E\cro{\| C^{1/2}(\tilde\theta-\theta)\|^2}
&\leq& c(K,\eta)\ \min_{m\in\M}\Biggl\{\| C^{1/2}(\theta-\theta_{m})\|^2\pa{1+{\pen({\rm  deg}(m))\over n-{\rm
deg}(m)}}\nonumber\\
&&{}+ {1\over n}\sum_{j=1}^p  (\pen(|m_{j}|)+K\log n)\sigma^2_{j} \Biggr\}+R_{n}(\eta,C)
\end{eqnarray}
where   $K$ is the constant appearing in~\eref{new},
$c(K,\eta)={K\over (K-1)(1-\sqrt\eta)^{4}}$ and the residual term $R_{n}(\eta,C)$ (made explicit in the proof) is of order a $p^2n^{-4\log n}$.
\end{theorem}

The proof of Theorem~\ref{main} and of the next corollary is delayed to Section~\ref{proofmain}.

\begin{corollary}\label{thm1}
Assume that $p>n\geq 3$ and that Condition~\eref{conditionD} holds. Then, there exists some constant $C_{K,\eta}$, depending on $K$ and $\eta$ only, such that
\begin{equation} \label{oraclethm1}
{\rm MSEP}(\tilde \theta)\leq C_{K,\eta}\log(p)\times\pa{ \min_{m\in\M}\acb{{\rm MSEP}(\thetah_{m})} \vee {1\over n}{\|C^{1/2}(I-\theta)\|^2}}+R_{n},
\end{equation}
where  $R_{n}=R_{n}(\eta,C)$ is of order a $p^2n^{-4\log n}$.
\end{corollary}

Corollary~\ref{thm1} roughly states that when the candidate graphs have a degree smaller than  $n/(2\log p)$,  the MSEP of $\tilde \theta$ nearly achieves, up to a $\log (p)$ factor, the minimal MSEP of the collection of estimators $\{\thetah_{m},\ m\in\M\}$.  In particular, if $\g\in\M$, the MSEP of $\tilde \theta$ is  upper-bounded by $\log(p)$ times  the MSEP of $\thetah_{\g}$, which in turn is roughly upper bounded by ${\rm deg}(\g)\times\|C^{1/2}(I-\theta)\|^2\log(p)/n$.

The additional term $n^{-1}\|C^{1/2}(I-\theta)\|^2$ in~\eref{oraclethm1} can be interpreted as a minimal variance for the estimation of $\theta$. This minimal variance is due to the inability of the procedure to detect  with probability one whether an isolated vertex of $\g$ is isolated or not. We mention that when each vertex of the graph $\g$ is connected to at least one other vertex, this variance term $n^{-1}\|C^{1/2}(I-\theta)\|^2$ remains smaller than the MSEP of $\thetah_{\g}$.

 Below, we discuss on the necessity of Condition~\eref{conditionD}  on the degree of the graphs and on the size of the penalty.

\subsection[Is Condition (5) avoidable?]{Is Condition \eref{conditionD} avoidable?}\label{optimal}

Condition~\eref{conditionD} requires that $D_{\M}$ remains small compared to $n/(2\log p)$. We may wonder if this condition is necessary, or if we can hope to handle graphs with larger degree $D$. A glance at the proof of Theorem~\ref{main} shows that Condition~\eref{conditionD} can be replaced by the weaker condition $\pa{\sqrt{D_{\M}+1}+\sqrt{2\log C_{p-1}^{D_{\M}}+{1/(4C_{p-1}^{D_{\M}})}}}^2\leq{\eta n}$. Using the classical bound $C_{p-1}^D\leq (ep/D)^{D}$, we obtain that the latter condition is satisfied when
\begin{equation}\label{conditionDbis}
D_{\M}\leq {\eta\over 3}\times {n\over 2.1+\log{p\over D_{\M}}},
\end{equation}
so we can replace Condition~\eref{conditionD} by Condition~\eref{conditionDbis} in Theorem~\ref{main}. Let us check now that we cannot improve (up to a multiplicative constant) upon \eref{conditionDbis}.

Phythagorean equality gives
$\|C^{1/2}(\theta-\thetah)\|^2=\|C^{1/2}(I-\thetah)\|^2-\|C^{1/2}(I-\theta)\|^2$,
so there is no hope to control the size of
$\|C^{1/2}(\theta-\thetah)\|^2$ if we do not have for some
$\delta\in(0,1)$ the inequalities
\begin{eqnarray}\label{Dmatrix}
\hspace*{-20pt}
&&(1-\delta)\|C^{1/2}(I-A)\|_{p\times p}\nonumber\\
\hspace*{-20pt}&&\quad \leq  {1\over \sqrt{n}}\,\|X(I-A)\|_{n\times p}
\leq (1+\delta)\|C^{1/2}(I-A)\|_{p\times p}\quad \textrm{for all }
A\in\bigcup_{m\in\M}\Theta_{m}\quad \quad
\end{eqnarray}
with large probability. Under Condition~\eref{conditionD}
or~\eref{conditionDbis}, Lemma~\ref{concentration} Section~\ref{proofs}
ensures that these inequalities hold  for any $\delta>\sqrt{\eta}$ with
probability $1-2\exp(-n(\delta-\sqrt \eta)^2/2)$. We emphasize next
that in the simple case where $C=I$, there exists a constant
$c(\delta)>0$ (depending on $\delta$ only) such that the
Inequalities~\eref{Dmatrix} cannot hold if  $\M^{\#}_{D}\subset \M$ or
$\M^{\#,+}_{D}\subset \M$ with
$$D\geq c(\delta) {n\over 1+\log{p\over n}}.$$
Indeed, when $C=I$ and $\M^{\#}_{D}\subset \M$ (or
$\M^{\#,+}_{D}\subset \M$), the Inequalities~\eref{Dmatrix} enforces
that $n^{-1/2}X$ satisfies the so-called $\delta$-Restricted Isometry
Property of order $D$ introduced by Cand\`es and Tao~\cite{CT}, namely
$$(1-\delta)\|\beta\|_{p\times 1}\leq \|n^{-1/2}X\beta\|_{p\times p}\leq (1+\delta)\|\beta\|_{p\times 1}$$
for all $\beta$ in $\R^p$ with at most $D$ non-zero components.
Barabiuk {\it et al.}~\cite{BDDVW} (see also Cohen {\it et
al.}~\cite{CDDV}) have noticed that there exists some constant
$c(\delta)>0$ (depending on $\delta$ only) such that no $n\times p$
matrix can fulfill the $\delta$-Restricted Isometry Property of order
$D$ if $D\geq c(\delta)n/(1+\log(p/n))$. In particular, the matrix $X$
cannot satisfies the Inequalities~\eref{Dmatrix} when
$\M^{\#}_{D}\subset \M$ (or $\M^{\#,+}_{D}\subset \M$)  with $D\geq
c(\delta)n/(1+\log(p/n))$.

\subsection{Can we choose a smaller penalty?}\label{min}

As mentioned before, under Condition~\eref{conditionD} the penalty
$\pen(d)$ given by~\eref{new} is approximately upper bounded by  $K
\pa{1+e^{\eta}\sqrt{2\log p}}^2(d+1)$. Similarly to Theorem~1 in Baraud
{\it et al.}~\cite{BGH}, a slight variation of the proof of
Theorem~\ref{main} enables to justify  the use of a penalty of the form
$\pen(d)=2Kd\log(p-1)$ with $K>1$ as long as $D_{\M}$ remains small
(the condition on $D_{\M}$ is then much stronger than
Condition~\eref{conditionD}). We underline in this section, that it is
not recommended to choose a smaller penalty. Indeed, next proposition
shows that  choosing a penalty of the form $\pen(d)=2(1-\gamma)d\log
(p-1)$ for some $\gamma\in(0,1)$ leads to a strong overfitting in the
simple case where  $\theta=0$, which corresponds to  $C=I$.

\begin{proposition}\label{minimal}
Consider three integers $1\leq D<n<p$ such that $p\geq
e^{2/(1-\gamma)}+1$ and $\M^{\#}_{D}\subset \M$ or
$\M^{\#,+}_{D}\subset \M$. Assume that   $\pen(d)=2(1-\gamma)d\log
(p-1)$ for some $\gamma\in(0,1)$  and $\theta=0$. Then, there exists
some constant $c(\gamma)$ made explicit in the proof, such that when
$\hat m$ is selected according to~\eref{selection}
$$\p\pa{|\hat m|\geq  {{c(\gamma)\min(n,p^{\gamma/4})\over (\log p)^{3/2}}\wedge \lfloor\gamma D/ 8\rfloor}}\geq 1-3(p-1)^{-1}-2e^{-\gamma^2n/8^3}.$$

In addition, in the case where  $\M=\M^{{\rm deg},+}_{D}$,  we have
\begin{eqnarray*}
&&\p\pa{|\hat m_{j}|\geq  {{c(\gamma)\min(n,p^{\gamma/4})\over (\log p)^{3/2}}\wedge \lfloor\gamma D/ 8\rfloor}}\\
&&\quad \geq 1-3(p-1)^{-1}-2e^{-\gamma^2n/8^3}\quad \textrm{ for all $j\in\ac{1,\ldots,p}$.}
\end{eqnarray*}
\end{proposition}

\section{Numerical study}\label{numerique}

In this section, we carry out a small simulation study to evaluate the
performance of our procedure. Our study concerns the behaviour of the
estimator $\thetah$ when the sparsity decreases
(Section~\ref{sparsity}) or when the number of covariates $p$ increases
(Section~\ref{covariates}). In this direction, we fix the sample size
$n$ to 15 (a typical value in post-genomics) and run simulations for
different values of $p$ and for different sparsity levels. For
comparison, we include the procedure ``or'' of Meinshausen and
B\"uhlmann~\cite{MB}. This choice is based on the numerical study of
Villers {\it et al.}~\cite{HSV}, where this procedure achieves a good
trade-off between the power and the FDR. We write henceforth ``MB'' to
refer to this procedure.

\subsection{Simulation scheme}

The graphs $\g$ are sampled according to the Erd\"os-R\'enyi  model:
starting from a graph with $p$ vertices and no edges, we set edges
between each couple of vertices  at random with probability $q$
(independently of the others). Then, we associate with a graph $\g$ a
positive-definite matrix $K$ with shape given by  $\g$ as follows. For
each $(i,j)\in\g$, we  draw $K_{i,j}=K_{j,i}$ from the uniform
distribution in $[-1,1]$ and set the elements on the diagonal of $K$ in
such a way that $K$ is diagonal dominant, and thus positive definite.
Finally, we normalize $K$ to have ones on the diagonal and set
$C=K^{-1}$.

For each value of $p$ and $q$ we sample 20 graphs and covariance
matrices $C$. Then, for each covariance matrix $C$, we generate 200
independent samples $(X_{1},\ldots,X_{15})$ of size 15 with law
$\p_{C}$.  For each sample, we estimate $\theta$ with our procedure and
the procedure of Meinshausen and B\"uhlmann. For our procedure, we set
$\M=\M^{{\rm deg}}_{4}$ and $K=2$ or $2.5$. For Meinshausen and
B\"uhlmann's estimator $\hat\theta_{{\rm MB}}$ we set $\lambda$
according to (9) in~\cite{MB} with $\alpha=5\%$, as recommended by the
authors.

 On the basis of the 20*200 simulations we evaluate the risk ratio
 $${\rm r.Risk}={{\rm MSEP}(\thetah)\over \min_{m}{\rm MSEP}(\thetah_{m})},$$
as well as the power  and the FDR for the detection of the edges of the
graph $\g$. The calculations are made with R \href{http://www.r-project.org/}{www.r-project.org/}.

\begin{table}[b]
\vspace*{-6pt}
\renewcommand{\arraystretch}{1.2}
\caption{Our procedure with $K=2$, $K=2.5$ and MB procedure:
Risk ratio (r.Risk), Power and FDR when $n=15$, $p=10$ and $q=10\%$, $30\%$ and $33\%$. }
\label{t1}
\begin{tabular}{|c|ccc|ccc|ccc|}
\hline
& \multicolumn{3}{c|}{$q=10\%$} & \multicolumn{3}{c|}{$q=30\%$}&\multicolumn{3}{c|}{$q=33\%$}\\ \hline
Estimator & r.Risk  & Power & FDR & r.Risk & Power & FDR& r.Risk & Power & FDR\\ \hline
$K=2$ & 2.3 & 82\% &4.9\% & 4.3 & 23\% & 6.8\%&4.4 & 13\% & 5.6\% \\ \hline
$K=2.5$ &2.5& 81\%&4.4\%&4.9&20\%&5.4\%&4.9 & 10\% & 4.1\% \\ \hline
MB & 3.3&81\%&3.7\%&6.9& 14\%& 2.9\% & 6.4 & 3.8\% & 1.1\%\\ \hline
\end{tabular}
\end{table}

\begin{table}[b!]
\vspace*{-6pt}
\caption{Our procedure with $K=2$, $K=2.5$ and MB procedure:
Risk ratio (r.Risk), Power and FDR when $n=15$, $s=1$ and $p=15$, $20$ and $40$. }
\label{t2}
\begin{tabular}{|c|ccc|ccc|ccc|}
\hline
& \multicolumn{3}{c|}{$p=15$} & \multicolumn{3}{c|}{$p=20$}&\multicolumn{3}{c|}{$p=40$}\\ \hline
Estimator & r.Risk  & Power & FDR & r.Risk & Power & FDR& r.Risk & Power & FDR\\ \hline
$K=2$ & 3.6 & 74\% &6.6\% &3.7  & 69\% & 6\%&5.4 & 68\% &5.4 \% \\ \hline
$K=2.5$ &4.3& 72\%& 6\%&4.4 &68\%&5.3\%&6.5 &67\% & 4.7\% \\ \hline
MB & 17 &60\%&4\%&160 & 20\%& 4.8\% & 340 & 0.0\% & 0.0\%\\ \hline
\end{tabular}
\end{table}

\vspace*{-3pt}
\subsection{Decreasing the sparsity}\label{sparsity}

To investigate the behaviour of the procedure when the sparsity
decreases, we fix $(n,p)=(15,10)$ and consider the three graph-density
levels $q=10\%$, $q=30\%$ and $q=33\%$. The results are reported in
Table~\ref{t1}.

When $q=10\%$ the procedures have a good performance. They detect on
average more than $80\%$ of the edges with a FDR lower than $5\%$ and a
risk ratio around 2.5. We note that MB has a slightly larger risk ratio
than our procedure, but also a slightly smaller FDR.

When $q$ increases above $30\%$ the performances of the procedures
decline abruptly. They detect less than $25\%$ of the edges on average
and the risk ratio increases above 4. When $q=30\%$ or $q=33\%$ our
procedure is more powerful than MB, with a risk ratio $33\%$ smaller.

In this simulation study, all the candidate graphs have a degree
smaller than 4. Using candidate graphs with a larger degree should not
change the nature of the results. Actually, when $q=30$ or $33\%$, less
than $2\%$ of the selected graphs have a degree equal to 4 and the mean
degree of the selected graphs is between 1 and 2.

\vspace*{-3pt}
\subsection{Increasing the number of covariates}\label{covariates}

In this section, we focus on the quality of the estimation of $\theta$
and $\g$ when the number of covariates $p$ increases. We thus fix the
sample size $n$ to 15 and the sparsity index $s:=pq$ to 1. This last
index corresponds to the mean degree of a vertex in the Erd\"os-R\'enyi
model. Then, we run simulations for three values of $p$, namely $p=15$,
$p=20$ and $p=40$ (in this last case we set $\M=\M^{{\rm deg}}_{3}$ to
reduce the computational time). The results are reported in
Table~\ref{t2}.

When the number $p$ of covariates increases, the risk ratios of the
procedures increase and their power decrease. Nevertheless, the
performance of our procedure remains good, with a risk ratio  between
3.6 and 6.5, a power close to $70\%$ and a FDR around $5.6\pm1\%$. In
contrast, the performances of MB decrease abruptly when $p$ increases.
For values of $p$ larger or equal to 22 (not shown), MB procedure does
not detect any edge anymore. This phenomenon was already noticed in
Villers {\it et al.}~\cite{HSV}.

\section{Conclusion}\label{conclusion}

In this paper, we propose to estimate the matrix of regression
coefficients $\theta$ by minimizing some penalized empirical risk. The
resulting estimator has some nice theoretical and practical properties.
From a theoretical point of view, Theorem~\ref{thm1} ensures that the
MSEP of the estimator can be upper-bounded in terms of the minimum of
the MSEP of the $\{\thetah_{m},\ m\in\M\}$ in a non-asymptotic setting
and with no condition on the covariance matrix $C$. From a more
practical point of view, the simulations of the previous section
exhibit a good behaviour of the estimator. The power and the risk of
our procedure are better than those of the procedure of Meinshausen and
B\"uhlmann, especially when $p$ increases. The downside of this better
power is a slightly higher FDR of our procedure compared to that of
Meinshausen and B\"uhlmann. If the FDR should be reduced, we recommend
to set the tuning parameter $K$ to a larger value, e.g.~$K=3$.

The main drawback of our procedure is its computational cost and in
practice it cannot be used when $p$ is larger than 50. In a future
work~\cite{GHV}, we propose a modification of the procedure that
enables to handle much larger values of $p$.

Finally, we emphasize that our procedure can only estimate accurately
graphs with a degree smaller than $n/(2\log p)$ and as explained in
Section~\ref{optimal}, we cannot improve (up to a constant) on this
condition.

\section{Proofs}\label{proofs}

\subsection{A concentration inequality}
\begin{lemma}\label{concentration}
Consider three integers $1\leq d\leq n\leq p$, a collection
$V_{1},\ldots,V_{N}$ of $d$-dimensional linear subspaces  of $\R^p$
and a $n\times p$ matrix $Z$ whose coefficients are i.i.d. with
standard gaussian distribution. We set
$\|\cdot\|_{n}=\|\cdot\|_{n\times 1}/\sqrt{n}$ and
$$\lambda^*_{d}(Z)=\inf_{v\in V_{1}\cup\cdots\cup V_{N}} {\|Z v\|_{n} \over  \ \|v\|_{p\times 1}}.$$
Then, for any $x\geq 0$
\begin{equation}\label{eq-concentration}
\p\pa{\lambda^*_{d}(Z)\leq 1-{\sqrt{d}+\sqrt{2\log N}+\delta_{N}+x\over \sqrt{n}}}\leq \p\pa{\NN\geq x}\leq e^{-x^2/2},
\end{equation}
where $\NN$ has a standard Gaussian distribution and $\delta_{N}=\pa{N\sqrt{8\log N}}^{-1}$.

Similarly, for any $x\geq 0$
\begin{equation}\label{eq-concentration2}
\p\pa{ \sup_{v\in V_{1}\cup\cdots\cup V_{N}} {\|Z v\|_{n} \over \ \|v\|_{p\times 1}}
\,{\geq}\, 1\,{+}\,{\sqrt{d}\,{+}\,\sqrt{2\log N}\,{+}\,\delta_{N}\,{+}\,x\over \sqrt{n}}}\leq \p\pa{\NN\geq x}\leq e^{-x^2/2}.
\end{equation}
\end{lemma}

\begin{proof}
The map $Z\to (\sqrt{n}\,\lambda^*_{d}(Z))$ is $1$-Lipschitz, therefore the Gaussian concentration inequality enforces that
$$\p\pa{\lambda^*_{d}(Z)\leq \E\pa{\lambda^*_{d}(Z)}-x/\sqrt{n}}\leq \p\pa{\NN\geq x}\leq e^{-x^2/2}.$$
To get \eref{eq-concentration}, we need to bound $ \E\pa{\lambda^*_{d}(Z)}$ from below.
For $i=1,\ldots,N$, we set
$$\lambda_{i}(Z)=\inf_{v\in V_{i}} {\|Z v\|_{n} \over \|v\|}.$$
We get from \cite{DS} the bound
$$\p\pa{\lambda_{i}(Z)\leq 1-\sqrt{d\over n}-{x\over \sqrt n}}\leq \p(\NN\geq x),$$
hence there exists some standard Gaussian random variables $\NN_{i}$ such that
$$\lambda_{i}(Z)\geq 1-\sqrt{d/n}-\pa{\NN_{i}}_{+}/\sqrt{n},$$
where $(x)_{+}$ denotes the positive part of $x$.
Starting from Jensen's inequality, we have for any $\lambda>0$
\begin{eqnarray*}
\E\paB{\max_{i=1,\ldots,N}(\NN_{i})_{+}}&\leq&{1\over \lambda}\,\log \E\pab{e^{\lambda\max_{i=1,\ldots,N}(\NN_{i})_{+}}}\\
&\leq& {1\over \lambda}\,\log\pa{\sum_{i=1}^N \E\pa{e^{\lambda(\NN_{i})_{+}}}}\\
&\leq& {1\over \lambda}\,\log N+ {1\over \lambda}\,\log\pab{e^{\lambda^2/2}+1/2}\\
&\leq& {\log N\over \lambda}+{\lambda\over 2}+{e^{-\lambda^2/2}\over2\lambda}.
\end{eqnarray*}
Setting $\lambda=\sqrt{2\log N}$, we finally get
\begin{eqnarray*}
 \E\pa{\lambda^*_{d}(Z)}&=&\E\paB{\min_{i=1,\ldots,N}\lambda_{i}(Z)}\ \geq\ 1-{\sqrt{d}+\sqrt{2\log N}+\delta_{N}\over \sqrt{n}}
\end{eqnarray*}
This concludes the proof of~\eref{eq-concentration} and the proof of~\eref{eq-concentration2} is similar.
\end{proof}

\subsection[Proof of Corollary 1]{Proof of Corollary \ref{thm1}}\label{proofthm1}\
Corollary~\ref{thm1} is a direct consequence of Theorem \ref{main} and of the three following facts.
\begin{enumerate}
\item The equality  $\sum_{j=1}^p\sigma_{j}^2=\| C^{1/2}(I-\theta)\|^2$
holds. \item Proposition 4 in Baraud {\it et al.}~\cite{BGH} ensures
that when $D_{\M}$ fulfills  Condition~\eref{conditionD}, there exists
a constant $C(K,\eta)$ depending on $K$ and $\eta$ only, such that
$${\pen(d) \over n-d}\leq C(K,\eta) \quad \textrm{for all } d\leq D_{\M}.$$

\item When $D_{\M}$ fulfills~\eref{conditionD} the MSEP of the estimator $\thetah_{m}$ is bounded from below by
$$\E\pa{\| C^{1/2}(\theta-\thetah_{m})\|^2}\geq \| C^{1/2}(\theta-\theta_{m})\|^2
+{1\over \pab{1+\sqrt{\eta/(2\log p)}}^2}\ \sum_{j=1}^p|m_{j}|\,{\sigma_{j}^2\over n}.$$
The latter inequality follows directly from Lemma~\ref{concentration}.
\end{enumerate}

Finally, to give an idea of the size of $C(K,\eta)$, we mention the following approximate bound (for $n$ and $p$ large)
$$C(K,\eta)={\pen(D_{\M})\over n-D_{\M}} \lesssim {K \pa{1+e^{\eta}\sqrt{2\log p}}^2 \over n-D_{\M}} \times  \eta \;{n\over 2\pa{1.1+\sqrt{\log p}}^2}\asymp K\eta \,e^{2\eta}.$$

\subsection[Proof of Theorem 1]{Proof of Theorem \ref{main}}\label{proofmain}\
The proof is split into two parts.

First, we bound from above $\E\crob{\|
C^{1/2}(\tilde\theta-\theta)\|^2}$ by $\pa{1-\sqrt \eta}^{-4}\E\crob{\|
X(\thetah-\theta)\|_{n}^2}+R_{n}$. Then, we bound this last term  by
the right hand side of~\eref{oracle1}.

To keep formulas short, we write henceforth $D$ for $D_{\M}$.

\noindent {\bf a. From $\E\crob{\| C^{1/2}(\tilde\theta-\theta)\|^2}$ to $\E\crob{\| X(\thetah-\theta)\|_{n}^2}$.}

We set $\|\cdot\|_{n}=\|\cdot\|_{n\times 1}/\sqrt{n}$, $\ \lambda_{0}=\pa{1-\sqrt \eta}^2$,
$$\lambda^1_{j}={\|X\theta^{(j)}\|_{n}\over \|C^{1/2}\theta^{(j)}\|}\quad \textrm{and}\quad  \lambda^*_{j}= \inf\ac{{\|XC^{-1/2}v\|_{n}\over \|v\|}:\ {v\in \bigcup_{m\in\M_{j,D}^*}V_{m}}} $$
where $V_{m}=C^{1/2}<\theta^{(j)}>+C^{1/2}\Theta^{(j)}_{m}$ and
$\M^*_{j,D}$ is the set of those subsets $m$ of
$\ac{1,\ldots,j-1,j+1,\ldots,p}\times\ac{j}$ with cardinality $D$.
Then, for any $j=1,\ldots,p$
\begin{eqnarray*}
\E\cro{\| C^{1/2}(\tilde\theta^{(j)}-\theta^{(j)})\|^2}&=& \E\cro{\| C^{1/2}(\thetah^{(j)}-\theta^{(j)})\|^2{\bf 1}_{\ac{\lambda^*_{j}\geq \lambda_{0},\ \thetah^{(j)}=\tilde\theta^{(j)}}}}\\
&&{}+\E\cro{\| C^{1/2}\theta^{(j)}\|^2{\bf 1}_{\ac{\lambda^*_{j}\geq \lambda_{0},\ \tilde\theta^{(j)}=0,\ \lambda^1_{j}\leq 3/2}}}\\
&&{}+\E\cro{\| C^{1/2}\theta^{(j)}\|^2{\bf 1}_{\ac{\lambda^*_{j}\geq \lambda_{0},\ \tilde\theta^{(j)}=0,\ \lambda^1_{j}> 3/2}}}\\
&&{}+\E\cro{\| C^{1/2}(\tilde\theta^{(j)}-\theta^{(j)})\|^2{\bf 1}_{\ac{\lambda^*_{j}< \lambda_{0}}}}\\
&=& \E_{1}^{(j)}+\E_{2}^{(j)}+\E_{3}^{(j)}+\E_{4}^{(j)}.
\end{eqnarray*}
We prove in the next paragraphs that $\sum_{j=1}^p\E_{1}^{(j)}\leq
\lambda_{0}^{-2}\E\crob{\| X(\thetah-\theta)\|_{n}^2}$ and that the
residual term
$R_{n}(\eta,C)=\sum_{j=1}^p(\E_{2}^{(j)}+\E_{3}^{(j)}+\E_{4}^{(j)})$ is
of order  $p^2T_{n}^{-2}$. The proofs of these bounds bear the same
flavor as the proof of Theorem~1 in Baraud~\cite{B}.

\noindent{\it Upper bound on $\E_{1}^{(j)}$.} Since
$$C^{1/2}(\thetah^{(j)}-\theta^{(j)})\in \bigcup_{m\in\M_{j,D}^*}V_{m},$$
 we have
$$\| C^{1/2}(\thetah^{(j)}-\theta^{(j)})\|^2{\bf 1}_{\ac{\lambda^*_{j}\geq \lambda_{0}}}\leq \lambda_{0}^{-2}\,\| X(\thetah^{(j)}-\theta^{(j)})\|_{n}^2$$
and therefore
\begin{equation}\label{E1}
\E_{1}^{(j)} \leq \lambda_{0}^{-2}\,\E\cro{\| X(\thetah^{(j)}-\theta^{(j)})\|_{n}^2}.
\end{equation}

\noindent{\it Upper bound on $\E_{2}^{(j)}$.} All we need is to bound
$\p\pab{{\lambda^*_{j}\geq \lambda_{0},\ \tilde\theta^{(j)}=0,\
\lambda^1_{j}\leq 3/2}}$ from above. Writing $\lambda^-$ for the
smallest eigenvalue of $C$, we have on the event $\ac{\lambda^*_{j}\geq
\lambda_{0}}$
$$\|\thetah^{(j)}\|\leq {\|C^{1/2}\thetah^{(j)}\|\over \sqrt{\lambda^-}}\leq {\|X\thetah^{(j)}\|_{n}\over \lambda_{0}\sqrt{\lambda^-}}.$$
Besides, for any $m\in\M$,
$$X\thetah^{(j)}_{m}=\textrm{Proj}_{X\Theta^{(j)}_{m}}\pab{X\theta^{(j)}+\sigma_{j}\eps^{(j)}}$$
with $\eps^{(j)}$ distributed as a standard Gaussian random variable in $\R^n$.

Therefore, on the event $\acb{\lambda^*_{j}\geq \lambda_{0},\ \tilde\theta^{(j)}=0,\ \lambda^1_{j}\leq 3/2}$ we have
\begin{eqnarray*}
\|\thetah^{(j)}\|&\leq& {  \|X\theta^{(j)}\|_{n}+\sigma_{j}\|\eps^{(j)}\|_{n}  \over \lambda_{0}\sqrt{\lambda^-}  }\\
&\leq& { 1.5\, \|C^{1/2}\theta^{(j)}\|+\sigma_{j}\|\eps^{(j)}\|_{n}  \over \lambda_{0}\sqrt{\lambda^-}  }.
\end{eqnarray*}
As a consequence,
\begin{eqnarray*}
\lefteqn{\p\pa{\lambda^*_{j}\geq \lambda_{0},\ \tilde\theta^{(j)}=0,\ \lambda^1_{j}\leq 3/2}}\\ &\leq& \p\pa{ { 1.5\, \|C^{1/2}\theta^{(j)}\|+\sigma_{j}\|\eps^{(j)}\|_{n}  \over \lambda_{0}\sqrt{\lambda^-}  }>T_{n}\sqrt{p} }\\
&\leq& \left\{ \begin{array}{ll}
1 & \textrm{when  }3\, \|C^{1/2}\theta^{(j)}\|  > \lambda_{0}\sqrt{p\lambda^-}\,T_{n} \\
\p\pa{{ 2\sigma_{j}\|\eps^{(j)}\|_{n}  > \lambda_{0}\sqrt{p\lambda^-}  }\,T_{n}} & \textrm{else,}
\end{array}\right. \\
&\leq& \left\{ \begin{array}{ll}
{9\, \|C^{1/2}\theta^{(j)}\|^2  / (\lambda_{0}^2{\lambda^-}\,pT_{n}^2) } & \textrm{when  }3\, \|C^{1/2}\theta^{(j)}\|  > \lambda_{0}\sqrt{p\lambda^-}\,T_{n} \\
{ 4\sigma_{j}^2/( \lambda_{0}^2\lambda^-\,p  T_{n}^2)} & \textrm{else.}
\end{array}\right. \\
\end{eqnarray*}
Finally,
\begin{equation}\label{E2}
\E_{2}^{(j)}\leq \|C^{1/2}\theta^{(j)}\|^2\;{9\, \|C^{1/2}\theta^{(j)}\|^2  + 4\sigma_{j}^2\over \lambda_{0}^2{\lambda^-}\,pT_{n}^2 }.
\end{equation}

\noindent{\it Upper bound on $\E_{3}^{(j)}$.} We note that
$n\pa{\lambda^1_{j}}^2$ follows a $\chi^2$ distribution, with $n$
degrees of freedom. Markov inequality then yields the bound
$$\p\pa{\lambda^1_{j}>3/2}\leq\exp\pa{-\,{n\over 2}\pa{9/4-1-\log(9/4)}}\leq \exp(-n/5).$$
As a consequence, we have
\begin{equation}\label{E3}
\E_{3}^{(j)}\leq \|C^{1/2}\theta^{(j)}\|^2 \exp(-n/5).
\end{equation}

\noindent{\it Upper bound on $\E_{4}^{(j)}$.}
Writing $\lambda^+$ for the largest eigenvalue of the covariance matrix $C$, we  have
\begin{eqnarray*}
\E_{4}^{(j)}&\leq& 2\E\cro{\pa{\|C^{1/2}\theta^{(j)}\|^2+\|C^{1/2}\thetah^{(j)}\|^2}{\bf 1}_{\ac{\lambda^*_{j}< \lambda_{0}}}}\\
&\leq&2\pa{\|C^{1/2}\theta^{(j)}\|^2+\lambda^+ pT_{n}^2}\p\pa{\lambda^*_{j}< \lambda_{0}}.
\end{eqnarray*}
The random variable $Z=XC^{-1/2}$ is $n\times p$ matrix whose
coefficients are i.i.d.\;and have the standard Gaussian distribution.
The condition \eref{conditionD} enforces the bound
$${\sqrt{D+1}+\sqrt{2\log |\M^*_{j,D}|}+\delta_{ |\M^*_{j,D}|}\over \sqrt n}\leq \sqrt{\eta},$$
so Lemma \ref{concentration} ensures that
$$\p\pa{\lambda^*_{j}< \lambda_{0}}\leq \exp\pa{-n(1-\sqrt{\eta})\eta/2}$$
and finally
\begin{equation}\label{E4}
\E_{4}^{(j)}\leq 2\pa{\|C^{1/2}\theta^{(j)}\|^2+\lambda^+ pT_{n}^2}\exp\pa{-n(1-\sqrt{\eta})\eta/2}.
\end{equation}

\noindent{\it Conclusion.}
Putting together the bounds \eref{E1} to \eref{E4}, we obtain
\begin{equation}\label{lien}
\E\cro{\| C^{1/2}(\tilde\theta-\theta)\|^2}=\sum_{j=1}^p\E\cro{\| C^{1/2}(\tilde\theta-\theta)\|^2}\leq  \lambda_{0}^{-2}\,\E\cro{\| X(\thetah-\theta)\|_{n}^2}+R_{n}(\eta,C)
\end{equation}
with $R_{n}(\eta,C)=\sum_{j=1}^p(\E_{2}^{(j)}+\E_{3}^{(j)}+\E_{4}^{(j)})$ of order a $p^2T_{n}^{-2}=p^2n^{-4\log n}$.

\noindent {\bf b. Upper bound on $\E\cro{\| X(\thetah-\theta)\|_{n}^2}$.}
Let $m^{*}$ be an arbitrary index in $\M$. Starting from the inequality
\begin{eqnarray*}
&&\sum_{j=1}^p\pa{\|X^{(j)}-X\thetah_{\hat m}^{(j)}\|^2\times \pa{1+{\pen(|\hat m_{j}|) \over n-|\hat m_{j}|}}}\\
&&\quad \leq \sum_{j=1}^p\pa{\|X^{(j)}-X\thetah_{m^*}^{(j)}\|^2\times \pa{1+{\pen(|m^*_{j}|)\over n-|m^*_{j}|}}}
\end{eqnarray*}
 and following the same lines as in the proof of Theorem 2 in Baraud {\it et al.}~\cite{BGH}
we obtain for any $K>1$
\begin{eqnarray*}
&&{K-1\over K} \sum_{j=1}^p\| X(\thetah^{(j)}-\theta^{(j)})\|_{n}^2\\
&&\quad \leq
 \sum_{j=1}^p\cro{\| X(\theta^{(j)}-\bar \theta^{(j)}_{m^*})\|_{n}^2 + R^{(j)}_{m^{*}}+{\sigma^2_{j}\over n}\pa{K U^{(j)}_{\hat m_{j}}- \pen(|\hat m_{j}|){ V^{(j)}_{\hat m_{j}}  \over n-|\hat m_{j}|}      }  },
\end{eqnarray*}
where for any $m\in\M$ and $j\in\ac{1,\ldots,p}$
\begin{eqnarray*}
&&X\bar \theta^{(j)}_{m}=\textrm{Proj}_{X\Theta^{(j)}_{m}}(X\theta^{(j)}),\quad\\
&&\E\pa{R^{(j)}_{m}\,\big|\,X^{(k)}, k\neq j}\leq \pen(|m_{j}|)\cro{{\| X(\theta^{(j)}-\bar \theta^{(j)}_{m})\|_{n}^2
 \over n-|m_{j}|}+{\sigma^2_{j}\over n}}\textrm{ a.s.}
\end{eqnarray*}
and the two random variables $U^{(j)}_{m_{j}}$ and $V^{(j)}_{m_{j}}$
are  independent with a $\chi^2(|m_{j}|+1)$ and a $\chi^2(n-|m_{j}|-1)$
distribution respectively.  Combining this bound with Lemma~6 in Baraud
{\it et al.}~\cite{BGH}, we get
\begin{eqnarray*}
&&{K-1\over K} \E\crob{\| X(\thetah-\theta)\|_{n}^2}\\
&&\quad \leq
 \E\cro{\| X(\theta-\bar \theta_{m^*})\|_{n}^2}
 +\sum_{j=1}^p  \pen(|m^*_{j}|)\cro{{\E\crob{\| X(\theta^{(j)}-\bar \theta^{(j)}_{m^*})\|_{n}^2}
 \over n-|m_{j}^*|}+{\sigma^2_{j}\over n}}\\
&&\qquad {}  +K\sum_{j=1}^p{\sigma^2_{j}\over n} \sum_{m_{j}\in\M_{j}}(|m_{j}|+1)\\
&&\qquad {}\times \H\pa{|m_{j}|+1,n-|m_{j}|-1,{(n-|m_{j}|-1)\pen(|m_{j}|)\over K(n-|m_{j}|)}},
\end{eqnarray*}
 where $\M_{j}=\ac{m_{j},\ m\in\M}$.
 The choice \eref{new} of the penalty ensures that the last term is upper bounded by $K\sum_{j=1}^p\sigma_{j}^2\log(n)/n$.
We also note that  $\| X(\theta^{(j)}-\bar \theta^{(j)}_{m^*})\|_{n}^2\leq \| X(\theta^{(j)}- \theta^{(j)}_{m^*})\|_{n}^2$ for all $j\in\ac{1,\ldots,p}$ since $X\bar \theta^{(j)}_{m^*}=\textrm{Proj}_{X\Theta^{(j)}_{m^*}}(X\theta^{(j)})$. Combining this inequality with
$\E\cro{\| X(\theta^{(j)}- \theta^{(j)}_{m^*})\|_{n}^2}={\| C^{1/2}(\theta^{(j)}-\theta^{(j)}_{m^*})\|^2}$, we obtain
\begin{eqnarray}
\lefteqn{{K-1\over K} \E\cro{\| X(\thetah-\theta)\|_{n}^2}}\nonumber\\
&\leq&
 {\| C^{1/2}(\theta- \theta_{m^*})\|^2}
 +\sum_{j=1}^p  \pen(|m^*_{j}|)\cro{{{\| C^{1/2}(\theta^{(j)}- \theta^{(j)}_{m^*})\|^2} \over n-|m_{j}^*|}+{\sigma^2_{j}\over
 n}}\nonumber\\
&&{}  +K\sum_{j=1}^p{\sigma^2_{j}\over n}\,\log n \nonumber \\
 &\leq& {\| C^{1/2}(\theta- \theta_{m^*})\|^2}\pa{1+{\pen(D)\over n-D}}+\sum_{j=1}^p(\pen(|m_{j}|)+K\log n)
 {\sigma_{j}^2\over n}\label{Xoracle}\quad \qquad
 \end{eqnarray}

\noindent{\bf c. Conclusion.} The bound \eref{Xoracle} is true for any $m^*$, so combined with \eref{lien} it  gives
\eref{oracle1}.

\subsection[Proof of Proposition 1]{Proof of Proposition \ref{minimal}}
The proof of Proposition \ref{minimal} is based on the following Lemma.

Let us consider a $n\times p$ random matrix $Z$
whose coefficients $Z_{i}^{(j)}$ are i.i.d.\;with standard Gaussian
distribution and  a random variable $\eps$ independant of $Z$, with
standard Gaussian law in $\R^n$.

To any subset $s$ of $\ac{1,\ldots,p}$ we associate the linear space
$V_{s}={\rm span}\{e_{j},\ j\in s\}\subset \R^p$, where
$\ac{e_{1},\ldots,e_{p}}$ is the canonical basis of $\R^p$. We write
$Z\thetah_{s}={\rm Proj}_{ZV_{s}}(\eps)$, we denote by $\hat s_{d}$
the set of cardinality $d$ such that
\begin{equation}\label{shatd}
\|Z\thetah_{\hat s_{d}}\|^2=\max_{|s|=d}\|Z\thetah_{s}\|^2.
\end{equation}
and we define
$$\crit'(s)=\|\eps-Z\thetah_{s}\|^2\pa{1+{\pen(|s|)\over n-|s|}}.$$

\begin{lemma}\label{stylise}
Assume that $p\geq e^{2/(1-\gamma)}$ and $\pen(d)=2(1-\gamma)d\log p$.
We write $D_{n,p}$ for the largest integer smaller than
$$5D/6, \quad  {p^{\gamma/4}\over (4\log p)^{3/2}}\quad \textrm{and} \quad  {\gamma^2 n\over 512(1.1+\sqrt{\log p})^2}.$$
Then, the probability to have
$$\crit'(s)>\crit'(\hat s_{D_{n,p}}) \ \textrm{for all $s$ with cardinality smaller than }\gamma D_{n,p}/6$$ is  bounded from below by $1-3 p^{-1}-2\exp(-n\gamma^2/512).$
\end{lemma}
The proof of this lemma is technical and in a first time we only give a
sketch of it. For the details, we refer  to
Section~\ref{proof-stylise}.

\vspace*{3pt}
\noindent {\bf Sketch of the proof of Lemma \ref{stylise}.} We have
\begin{eqnarray*}
\|Z\thetah_{s}\|^2&=& \|\eps\|^2-\inf_{\hat \alpha\in V_{s}}\|\eps-Z\hat \alpha\|^2\\
&=& \sup_{\hat \alpha\in V_{s}}\cro{2<\eps,Z\hat \alpha>-\|Z\hat \alpha\|^2}.
\end{eqnarray*}
According to Lemma~\ref{concentration}, when $|s|$ is small compared to
$n/\log p$, we have $\|Z\hat \alpha\|^2\approx n \|\hat \alpha\|^2$
with large probability and then
$$\|Z\thetah_{s}\|^2\approx  \sup_{\hat \alpha\in V_{s}}\cro{2<Z^T\eps,\hat \alpha>-n\|\hat \alpha\|^2}={1\over n}\;\|{\rm Proj}_{V_{s}}(Z^T\eps)\|^2.$$
Now, $Z^T\eps=\|\eps\| Y$ with $Y$ independent of $\eps$ and with $\NN(0,I_{p})$ distribution, so
 $$\|Z\thetah_{s}\|^2\approx  {\|\eps\|^2\over n}\;\|{\rm Proj}_{V_{s}}Y\|^2.$$
Since $\max_{|s|=d}\|{\rm Proj}_{V_{s}}Y\|^2\approx 2d\log p$ with large probability, we have
$\|Z\thetah_{\hat s_{d}}\|^2\approx 2d\log p\times\|\eps\|^2/n$ and then
$$\min_{|s|=d}\crit'(s)=\crit'(\hat s_{d})\approx \|\eps\|^2\pa{1-{2\gamma d \log p \over n}}.$$
Therefore,  with large probability we have
$\crit'(s)>\crit'(\hat s_{D_{n,p}})$ for all $s$ with cardinality less than $\gamma D_{n,p}/6$.

\begin{proof}[Proof of Proposition \ref{minimal}] We start with the
case $\M^{\#,+}_{D}\subset \M$. When $|\hat m|\leq \gamma D_{n,p-1}/6$,
we have in particular $|\hat m_{1}|\leq \gamma D_{n,p-1}/6$. We build
$\tilde m$ from  $\hat m$ by replacing $\hat m_{1}$ by a set $\tilde
m_{1}\subset\ac{1}\times \ac{2,\ldots,p}$ which maximizes $\|X\hat
\theta_{\tilde m}^{(1)}\|^2$ among all the subset $\tilde m_{1}$ of
$\ac{1}\times \ac{2,\ldots,p}$ with cardinality $D_{n,p-1}$. It follows
from Lemma~\ref{stylise} (with $p$ replaced by $p-1$) that the
probability to have $\crit(\hat m)\leq \crit(\tilde m)$ is bounded from
above by $3 (p-1)^{-1}+2\exp(-n\gamma^2/512).$ Since $\tilde
m\in\M^{\#,+}_{D}$, the first part of Proposition~\ref{minimal}
follows. When $\M^{\#}_{D}\subset \M$, the proof is similar.

When $\M^{{\rm deg},+}_{D}\subset \M$, the same argument shows that for
any $j\in\ac{1,\ldots,p}$ the probability to have $|\hat m_{j}|\leq
\gamma D_{n,p-1}/6$ is bounded from above by $3
(p-1)^{-1}+2\exp(-n\gamma^2/512).$
\end{proof}

\subsection[Proof of Lemma 2]{Proof of Lemma \ref{stylise}}\label{proof-stylise}
We write $D$ for $D_{n,p}$ and  $\Omega_{0}$ for the event
$$\Omega_{0}=\left\{
\begin{array}{@{}ll@{}}
\|Z\thetah_{\hat s_{D}}\|^2\geq2D(1-\gamma/2)\|\eps\|_{n}^2\log p& \textrm{and}\\[3pt]
 \|Z\thetah_{s}\|^2\leq2\,|s|\,(2+\gamma)\|\eps\|_{n}^2\log p,& \textrm{for all } s \ \textrm{with } |s|\leq D
 \end{array}\right\}.$$
We will prove first that  on the event $\Omega_{0}$ we have
$\crit'(s)>\crit'(\hat s_{D_{n,p}})$ for any $s$ with cardinality less
than $\gamma D_{n,p}/6$ and then  we will prove that $\Omega_{0}$ has a
probability bounded from below by $1-3 p^{-1}-2\exp(-n\gamma^2/512)$.

We write $\Delta(s)=\crit'(\hat s_{D})-\crit'(s)$. Since we are
interested in the sign of $\Delta(s)$, we will still write $\Delta(s)$
for any positive constant times $\Delta(s)$. We have on~$\Omega_{0}$
\begin{eqnarray*}
{\Delta(s)\over\|\eps\|^2} &\leq& \pa{1-{2\log p\over n}(1-\gamma/2)D}\pa{1+{\pen(D)\over
n-D}}\\
&&{}-\pa{1-{2\log p\over n}(2+\gamma)| s|}\pa{1+{\pen(|s|)\over n-|s|}}.
\end{eqnarray*}
We note that ${{\pen(|s|)/( n-|s|)}}\leq {{\pen(D)/(n-D)}}$. Multiplying by $n/(2\log p)$ we obtain
\begin{eqnarray*}
\Delta(s)&\leq &(1-\gamma)D\pa{1+{D-2(1-\gamma/2)D\log p\over n-D}}-(1-\gamma/2)D\\
&&{}-(1-\gamma)|s|+(2+\gamma)|s|+(2+\gamma)|s|{\pen(D)\over n-D}\\
&\leq& (1-\gamma)D\pa{1+{D-2(1-\gamma/2)D\log p+2(2+\gamma)|s|\log p\over n-D}}\\
&&{}-(1-\gamma/2)D+(1+2\gamma)|s|.
\end{eqnarray*}
When $p \geq e^{2/(1-\gamma)}$ and $|s|\leq \gamma D/6$ the first term
on the right hand side is bounded from above by $(1-\gamma)D$, then
since $\gamma<1$
\begin{eqnarray*}
\Delta(s)&\leq & (1+2\gamma)\gamma D/6-\gamma D/2\ < \ 0.
\end{eqnarray*}

We will now bound $\p\pa{\Omega^c_{0}}$ from above.
We write $Y=Z^T\eps/\|\eps\|$ (with the convention that $Y=0$ when $\eps=0$) and
\begin{eqnarray*}
\Omega_{1}&=&\ac{{2\over 2+\gamma}\leq  {\|Z\hat \alpha\|^2_{n}\over\|\hat \alpha\|^2 }\leq \pa{1-\gamma/2}^{-1/2}, \quad \textrm{for all } \hat \alpha\in\bigcup_{|s|=D}V_{s}},\\
\Omega_{2}&=&\ac{ \max_{|s|=D}\|{\rm Proj}_{V_{s}}Y\|^2\geq 2(1-\gamma/2)^{1/2}D\log p   },\\
\Omega_{3}&=&\ac{\max_{i=1,\ldots,p}Y_{i}^2\leq 4\log p}.
\end{eqnarray*}
We first prove that $\Omega_{1}\cap\Omega_{2}\cap\Omega_{3}\subset \Omega_{0}$. Indeed, we have on $\Omega_{1}\cap \Omega_{2}$
\begin{eqnarray*}
\|Z\thetah_{\hat s_{D}}\|^2&=& \max_{|s|=D}\sup_{\hat \alpha\in V_{s}}\cro{2<\eps,Z\hat \alpha>-\|Z\hat \alpha\|^2}\\
&\geq&  \max_{|s|=D}\sup_{\hat \alpha\in V_{s}}\cro{2<Z^T\eps,\hat \alpha>-n(1-\gamma/2)^{-1/2}\|\hat \alpha\|^2}\\
&\geq&  {\pa{1-\gamma/2}^{1/2}\|\eps\|^2\over n}\max_{|s|=D}\|{\rm Proj}_{V_{s}}Y\|^2\\
&\geq& {2D\pa{1-\gamma/2}\|\eps\|^2_{n}\log p}.
\end{eqnarray*}
Similarly,  on $\Omega_{1}$ we have $\|Z\thetah_{s}\|^2\leq
\|\eps\|^2_{n}\|{\rm Proj}_{V_{s}}Y\|^2\times (2+\gamma)/2$ for all $s$
with cardinality less than $D$. Since $\|{\rm Proj}_{V_{s}}Y\|^2\leq
|s|\max_{i=1,\ldots,p}(Y_{i}^2)$, we have  on
$\Omega_{1}\cap\Omega_{3}$
$$\|Z\thetah_{s}\|^2\leq 2 (2+\gamma)|s|\,\|\eps\|^2_{n}\log p,$$
for all $s$ with cardinality less than $D$
and then $\Omega_{1}\cap\Omega_{2}\cap\Omega_{3}\subset \Omega_{0}$.

To conclude, we bound $\p(\Omega_{i}^c)$ from above, for $i=1,2,3$. First, we have
$$\p(\Omega_{3}^c)=\p\pa{\max_{i=1,\ldots,p}Y_{i}^2>4\log p}\leq 2p\,\p(Y_{1}\geq 2\sqrt{\log(p)})\leq 2p^{-1}.$$
To bound  $\p(\Omega_{1}^c)$, we note that $(1-\gamma/2)^{-1/4}\geq
1+\gamma/8$ and $\sqrt{2/(2+\gamma)}\leq 1-\gamma/8$ for any
$0<\gamma<1$, so Lemma~\ref{concentration} ensures that
$\p(\Omega_{1}^c)\leq 2e^{-n\gamma^2/512}$. Finally, to bound
$\p(\Omega_{2}^c)$, we sort the $Y_{i}^2$ in decreasing order
$Y_{(1)}^2>Y_{(2)}^2>\cdots>Y_{(p)}^2$ and note that
$$\max_{|s|=D}\|{\rm Proj}_{V_{s}}Y\|^2\geq DY^2_{(D)}.$$
Furthermore, we have
\begin{eqnarray*}
\p\pa{Y^2_{(D)}\leq 2(1-\gamma/2)^{1/2}\log p}&\leq &\binom{D-1}{p}\p\pa{Y_{1}^2\leq  2(1-\gamma/2)^{1/2}\log p} ^{p-D+1}\\
&\leq& p^{D-1}\pa{1-{p^{\sqrt{1-\gamma/2}}\over 4(1-\gamma/2)^{1/4}\sqrt{2\log p}}}^{p-D+1},
\end{eqnarray*}
where the last inequality follows from $p\geq e^{2/(1-\gamma)}$ and Inequality (60) in Baraud {\it et al.}~\cite{BGH}.
Finally, we obtain
\begin{eqnarray*}
\p\pa{Y^2_{(D)}\leq 2(1-\gamma/2)^{1/2}\log p}&\leq & p^{-1}\exp\pa{D\log p-{(p-D+1)p^{\sqrt{1-\gamma/2}}\over 4(1-\gamma/2)^{1/4}\sqrt{2\log p}}}\\
&\leq& p^{-1},
\end{eqnarray*}
where the last inequality comes from $ D\leq {p^{\gamma/4}/ (4\log p)^{3/2}}$.
To conclude $\p(\Omega_{2}^c)\leq p^{-1}$ and $\p\pa{\Omega_{0}^c}\leq
3 p^{-1}+2\exp(-n\gamma^2/512)$.

\end{document}